\documentclass[12pt]{amsart}
\usepackage{amsmath}
\usepackage{amscd}
\usepackage{diagrams}
\usepackage{amssymb}
\usepackage[mathscr,psamsfonts]{eucal}

\newcommand{\mysec}[1]{\section{#1}}

\newcommand{\myssec}[1]{\subsection{#1}}

\newtheorem{Caution}{Caution}[section]

\newtheorem{Corollary}[Caution]{Corollary}
\newtheorem{Definition}[Caution]{Definition}
\newtheorem{Example}[Caution]{Example}

\newtheorem{Lemma}[Caution]{Lemma}

\newtheorem{Note}[Caution]{Note}

\newtheorem{Proposition}[Caution]{Proposition}
\newtheorem{Remark}[Caution]{Remark}
\newtheorem{Theorem}[Caution]{Theorem}
\newcommand{\bCr}{\begin{Corollary}\em}
\newcommand{\eCr}{\end{Corollary}}
\newcommand{\bDf}{\begin{Definition}\em}
\newcommand{\eDf}{\end{Definition}}
\newcommand{\bEx}{\begin{Example}\em}
\newcommand{\eEx}{\end{Example}}
\newcommand{\bLm}{\begin{Lemma}\em}
\newcommand{\eLm}{\end{Lemma}}
\newcommand{\bNt}{\begin{Note}\em}
\newcommand{\eNt}{\end{Note}}
\newcommand{\bPf}{\begin{proof}[\noindent\indent{\sc Proof}]}
\newcommand{\ePf}{\end{proof}}
\newcommand{\bPr}{\begin{Proposition}\em}
\newcommand{\ePr}{\end{Proposition}}
\newcommand{\bRm}{\begin{Remark}\em}
\newcommand{\eRm}{\end{Remark}}
\newcommand{\bTh}{\begin{Theorem}}
\newcommand{\eTh}{\end{Theorem}}
\newcommand{\bEq}{\begin{eqnarray}}
\newcommand{\eEq}{\end{eqnarray}}
\newcommand{\beq}{\begin{eqnarray*}}
\newcommand{\eeq}{\end{eqnarray*}}
\newcommand{\bal}{\begin{align*}}
\newcommand{\bAl}{\begin{align}}
\newcommand{\bat}{\begin{alignat*}}
\newcommand{\fz}{\footnotesize}
\newcommand{\Rn}{{I\!\!R}}

\newcommand{\der}{\partial}

\newcommand{\Fla}{^{\flat}{}}

\newcommand{\nab}{\nabla}

\newcommand{\car}{\times}
\newcommand{\ten}{\otimes}

\newcommand{\wed}{\wedge}
\DeclareMathOperator{\con}{\lrcorner}

\DeclareMathOperator{\byd}{\,{\rm =}{\raisebox{.092ex}{\rm :}}\,}

\newcommand{\uten}[1]{\underset{#1}{\otimes}}

\newcommand{\fr}[2]{\frac{#1}{#2}\,}

\newcommand{\co}[2]{_{#1}{}^{#2}}
\newcommand{\col}[3]{_{#1}{}^{#2}{}_{#3}}

\newcommand{\cul}[4]{_{#1#2}{}^{#3}{}_{#4}}

\newcommand{\ENDE}{{\,\text{\footnotesize\qedsymbol}}}

\newcommand{\ssep}[1]{{\qquad\text{\rm{#1}}\qquad}}

\DeclareMathOperator{\map}{{{map}}}

\newcommand{\f}[1]{{\boldsymbol{#1}}}

\newcommand{\ba}[1]{{{\bar{#1}}}}

\newcommand{\dt}[1]{{\dot{#1}}}

\newcommand{\B}[1]{{\mathbb{#1}}}

\newcommand{\bet}{\beta}
\newcommand{\kap}{\kappa}
\newcommand{\lam}{\lambda}
\newcommand{\sig}{\sigma}
\newcommand{\ups}{\upsilon}
\newcommand{\Ups}{\Upsilon}
\newcommand{\ome}{\omega}
\newcommand{\Lam}{\Lambda}

\begin{document}
\title[Geometric structures]
{Geometric structures on the tangent bundle\\
        of the Einstein spacetime}

\author[J. Jany\v{s}ka]
{Josef Jany\v{s}ka}

\address{{\ }
\newline
Department of Mathematics, Masaryk University
\newline
Jan\'a\v{c}kovo n\'am 2a, 662 95 Brno, Czech Republic
\newline
\fz email: {\tt janyska@math.muni.cz}
}


\thanks{This research has been supported
by the Ministry of Education of the Czech Republic under the project
MSM0021622409 and
by the Grant agency of the Czech Republic under the project GA
201/05/0523.
}

\maketitle
\begin{abstract}
We describe conditions under which a spacetime connection
and a scaled Lorentzian metric define natural symplectic
and Poisson structures on the tangent bundle
of the Einstein spacetime.
\end{abstract}


\section*{Introduction}
\label{Introduction}

Geometrical structures induced on the tangent bundle of the Einstein
spacetime play a fundamental role in the covariant classical and quantum
mechanics.
The covariant classical and quantum  mechanics
over
the Einstein spacetime
proposed in
\cite{JanMod96a, JanMod97} is
natural in the sense of
\cite{KruJan90, KolMicSlo93, Nij72} and independent
on the base of scales, so  the ``spaces of scales''
are systematically used.
Roughly speaking, a space of scales
has the algebraic structure of
$\Rn^+$
but has no distinguished `basis'.
The basic objects of the theory (metric, 2-forms, 2-vectors, etc.)
are valued into {\em scaled\/} vector bundles,
that is into vector
bundles multiplied tensorially with spaces of scales.
In this way, each tensor field carries explicit information on its
``scale dimension".
Actually, in this paper, we assume the
space of {\em lengths}
$\B L \,.$
Moreover, $\B L^p$ denotes $\ten^p\B L$.

In \cite{Jan95, Jan00} the classification of symplectic and Poisson
structures on the tangent bundle of a pseudo-Riemannian
manifold was given for a non-scaled metric $g$ and a
torsion free linear connection $K$.
In this case the metric $g$ and the connection $K$
admit a family of symplectic 2-forms $\Ups[g,K]$
or Poisson 2-vectors $\Lam[g,K]$
on the tangent bundle parametrized by a function $\mu(g(u,u))$
satisfying certain conditions. Moreover, $g$ and $K$ are related
by the condition that $\nab g$ is a symmetric (0,3)-tensor field.
For a scaled metric $g$ and a general spacetime connections
the constructions of the 2-form $\Ups[g,K]$
and the 2-vector $\Lam[g,K]$ are the same but from the independence
on the base of scales it follows that $\Ups[g,K]$
and $\Lam[g,K]$ are
unique, up to a multiplicative real constant.
In this paper we generalize the results of \cite{Jan95, Jan00}
for a scaled metric and a general spacetime connection.
Namely, we shall describe
the condition under which the 2-form $\Ups[g,K]$
and the 2-vector $\Lam[g,K]$ give symplectic and
Poisson structures, respectively.

If
$\f E$
is a manifold, then
the tangent bundle will be denoted by $\tau[\f E]: T\f E\to \f E$
and local coordinates $(x^\lam)$ on $\f E$ induce the fibered
local coordinates $(x^\lam,\dt x^\lam)$ on $T\f E$.
By $\map (\f E,\f E')$ we denote the sheaf of smooth maps.

\smallskip


\mysec{Geometry of the spacetime}

We recall basec properties of the Einstein spacetime
and its tangent bundle.

\myssec{Spacetime}
\label{Spacetime}
We assume {\em spacetime\/} to be an oriented and time oriented
4--dimensional manifold
$\f E$
equipped with a scaled Lorentzian metric
$g : \f E \to \B L^2 \ten (T^*\f E \ten T^*\f E)$
with signature
$(-+++) \,.$
The dual metric will be denoted by
$\ba g : \f E \to \B L^{*2} \ten (T\f E \ten T\f E)$.
Let us note that the dimension is not relevant. Our results are
valid for any dimension $n\ge 3$ and a pseudo-Riemannian metric
of the signature $(1,n-1)$.

A {\em spacetime chart\/} is defined to be an ordered chart
$(x^0, x^i) \in \map(\f E, \, \Rn \car \Rn^3)$
of
$\f E \,,$
which fits the orientation of spacetime and such that the vector
$\der_0$
is timelike and time oriented and the vectors
$\der_1, \der_2, \der_3$
are spacelike.
In the following we shall always refer to spacetime charts.
Latin indices
$i, j, \dots$
will span spacelike coordinates, while Greek indices
$\lam, \mu, \dots$
will span spacetime coordinates.

We have the coordinate expressions
\bat{3}
g
&=
g_{\lam\mu} \, d^\lam \ten d^\mu \,,
&&\ssep{with}
g_{\lam\mu}
&&\in \map(\f E, \, \B L^2 \ten \Rn)
\\
\ba g
&=
g^{\lam\mu} \, \der_\lam \ten \der_\mu\,,
&&\ssep{with}
g^{\lam\mu}
&&\in \map(\f E, \, \B L^{*2}\ten\Rn) \,.
\end{alignat*}

\myssec{Spacetime connections}
\label{Spacetime connections}

We define a {\em (general) spacetime connection\/} to be a connection
$K$
of the bundle
$\tau[\f E] : T\f E \to \f E \,.$
We recall that a connection
$K$
of the bundle
$T\f E \to \f E$
can be expressed, equivalently, by a tangent valued form
$
K : T\f E \to T^*\f E \ten TT\f E \,,
$
which is projectable over
$\f 1 : \f E \to T^*\f E \ten T\f E \,,$
or by the vertical valued form
$
\nu[K]  :
T\f E \to T^*T\f E \ten VT\f E \,.
$
Their coordinate expressions are of the type
\bEq\label{Eq1.1}
K
=
d^\lam \ten (\der_\lam + K\co\lam\nu \, \dt\der_\nu)\,,
\qquad
\nu[K]
=
(\dt d^\nu - K\co\lam\nu \, d^\lam) \ten \dt\der_\nu \,,
\eEq
where
$K\co\lam\nu \in \map(T\f E, \, \Rn) \,$ and
$(\der_\lam,\dt \der_\lam)$ or $(d^\lam,\dt d^\lam)$
are the induced bases of local sections of
$TT\f E\to T\f E$ or $T^*T\f E\to T\f E$, respectively.

\smallskip

The connection
$K$
is said to be {\em linear\/} if it is a linear fibred morphism over
$\f 1 : \f E \to T^*\f E \ten T\f E \,.$
Moreover, the connection
$K$
is linear if and only if its
coordinate  expression is of the type
\beq
K\co\lam\nu = K\col\lam\nu\mu \, \dt x^\mu \,,
\ssep{with}
K\col\lam\nu\mu \in \map(\f E, \, \Rn) \,.
\eeq

The {\em torsion\/} of the connection
$K$
is defined to be the vertical valued 2--form
\beq
\tau[K] \byd - [\vartheta, K] : T\f E \to \Lam^2T^*\f E \ten VT\f E
\,,
\eeq
where
$[\,,]$
is the Fr\"olicher-Nijenhuis bracket and
$\vartheta : T\f E \to T^*\f E \ten VT\f E$
is the natural vertical valued 1--form with coordinate expression
$\vartheta = d^\lam \ten \dt\der_\lam \,.$
We have the coordinate expression
\bEq\label{Eq1.2}
\tau[K] =
\dt\der_\mu K\co \lam\nu \, d^\lam \wed d^\mu \ten \dt\der_\nu
\,.
\eEq

In the linear case, the torsion can be identified with a section
$
\tau[K] : \f E \to \Lam^2T^*\f E \ten T\f E
$
and its coordinate expression turns out to be the usual formula
$
\tau[K] = K\col \lam\nu\mu \, d^\lam \wed d^\mu \ten \der_\nu \,.
$
Thus, the connection
$K$
is linear and torsion free if and only if its coordinate expression
is of the type
\beq
K\co\lam\nu = K\col\lam\nu\mu \, \dt x^\mu \,,
\ssep{with}
K\col\lam\nu\mu = K\col\mu\nu\lam \in \map(\f E, \, \Rn) \,.
\eeq

We shall denote by $K[g]$ the canonical torsion free
linear spacetime metric connection given by $\nab g=0.$
We have
\bEq\label{Eq1.3}
    K[g]\col \mu\lam\nu =
        -\fr12 g^{\lam\rho}\,(\der_\mu g_{\rho\nu} +
        \der_\nu g_{\rho\mu} - \der_\rho g_{\mu\nu})\,.
\eEq

The {\em curvature\/} of the connection
$K$
is defined to be the vertical valued 2--form
\bEq\label{Eq1.4}
R[K] \byd - [K, K] : T\f E \to \Lam^2T^*\f E \ten VT\f E \,,
\eEq
where
$[\,,]$
is the Fr\"olicher-Nijenhuis bracket.
We have the coordinate expression
\bAl\label{Eq1.5}
R[K] & = R[K]_{\lam\mu}{}^\nu\,d^\lam \wed d^\mu \ten \dt\der_\nu
\\
      &  = - 2 \,
(\der_\lam K\co\mu\nu + K\co\lam\rho \, \dt\der_\rho K\co\mu\nu) \,
d^\lam \wed d^\mu \ten \dt\der_\nu
\,.        \nonumber
\end{align}

In the linear case, the coordinate expression turns out to be the
usual formula
\bAl\label{Eq1.6}
R[K] & = R[K]_{\lam\mu}{}^\nu{}_\sig \, \dt x^\sig
d^\lam \wed d^\mu \ten \dt\der_\nu
\\        \nonumber
&
= - 2 \,
(\der_\lam K\col\mu\nu\sig + K\col\lam\rho\sig \, K\col\mu\nu\rho) \,
\dt x^\sig
d^\lam \wed d^\mu \ten \dt\der_\nu \,.
\end{align}
Hence, in the linear case, the curvature can be identified with a
section
\beq
R[K] : \f E \to \Lam^2T^*\f E \ten T\f E \ten T^*\f E
\,,
\eeq
with the usual coordinate expression
\begin{align}
\label{Eq1.7}
R[K]
&=
R[K]\cul\lam\mu\nu\sig \, d^\lam \wed d^\mu
        \ten \der_\nu \ten d^\sig
\\ \nonumber
&= - 2 \,
(\der_\lam K\col\mu\nu\sig + K\col\lam\rho\sig \, K\col\mu\nu\rho) \,
d^\lam \wed d^\mu \ten \der_\nu \ten d^\sig \,.
\end{align}

\myssec{The Lie derivative and the exterior covariant differential
with respect to a spacetime connection}

A (general) spacetime connection $K$ considered as a tangent valued
1-form on $T\f E$ admits as usual, \cite{KolMicSlo93},
the Lie derivative of forms on
$T\f E$. Namely,
\beq
    L[K] \,\phi = \big( i(K) \,d - d\, i(K)\big)\, \phi:
        T\f E\to \Lam^{r+1}T^*T\f E
\eeq
for any $r$-form $\phi:T\f E\to \Lam^rT^*T\f E$.
Similarly we can define the Lie derivative
\beq
    L\big[R[K]\big] \,\phi = \big( i(R[K]) \,d +
        d\, i(R[K])\big)\, \phi:
        T\f E\to \Lam^{r+2}T^*T\f E\,.
\eeq

On the other hand a linear spacetime connection $K$
admits covariant exterior differential, \cite{KolMicSlo93},
of vector-valued forms on $\f E$. We apply this operation
on $T^*\f E$-valued forms on $\f E$ and compare
it with the Lie derivative.

Let $\phi$ be an $T^*\f E$-valued $r$-form on $\f E$,
or equivalently
$\phi:\f E\to \Lam^r T^*\f E\otimes_{\f E}T^*\f E$
be a section.
The {\it covariant exterior differential\/}
of $\phi$ with respect to $K$ is then
defined to be the $T^*\f E$-valued
$(r+1)$-form $d_K\phi$ on $\f E$ given by
\begin{gather} \label{Eq1.8}
d_K\phi(X_1,\cdots,X_{r+1})(Y)=\sum_{i=1}^{r+1}(-1)^{i+1}
\nabla_{X_i}(\phi (X_1,\cdots,\hat X_i,\cdots
,X_{r+1}))(Y)
\\            \nonumber
+\sum_{i<j}(-1)^{i+j}\phi([X_i,X_j],X_1,\dots,
\hat X_i,\dots,\hat X_j,
\dots,X_{r+1})(Y),
\end{gather}
for any vector fields $Y,X_1,\cdots,X_{r+1}$ on $\f E$.
The vector fields
$\hat X_i$ are omitted.

\smallskip

Any $T^*\f E$-valued $r$-form on $\f E$ can be considered to be a
linear horisontal $r$-form on $T\f E$. Then we have

\bLm\label{Lm1.1}
Let $\phi$ be a linear horisontal $r$-form on $T\f E$
and $K$ be a spacetime connection. Then
the Lie derivative $L[K]\, \phi$ is a linear horisontal
$(r+1)$-form on $T\f E$ if and only if $K$ is linear.
Moreover,
$L[K]\, \phi$ and $d_K \phi$ coincides.
\eLm

\bPf
Let $\phi$ is a linear horisontal $r$-form given in coordinates
by $\phi=\phi_{\rho\lam_1\dots\lam_r}\,\dt x^\rho\, d^{\lam_1}\wed
\dots \wed d^{\lam_r}$, $\phi_{\rho\lam_1\dots\lam_r}\in
\map(\f E,\Rn)$.
Then we have
\beq
    L[K]\, \phi =
        (\der_\mu \phi_{\rho\lam_1\dots\lam_r} \, \dt x^\rho
  + \phi_{\sig\lam_1\dots\lam_r}\, K_\mu{}^\sig{})\,
  d^\mu\wed d^{\lam_1} \wed\dots \wed d^{\lam_r}
        \,,
\eeq
i.e., in the linear spacetime connection case,
\beq
    L[K]\, \phi =
        (\der_\mu \phi_{\rho\lam_1\dots\lam_r}
  + \phi_{\sig\lam_1\dots\lam_r}\, K_\mu{}^\sig{}_\rho)\,\dt x^\rho\,
  d^\mu\wed d^{\lam_1} \wed\dots \wed d^{\lam_r}
        \,,
\eeq
which implies that $L[K]\, \phi$
is a linear horisontal $(r+1)$-form.

On the other hand $\phi$ can be considered to be a
$T^*\f E$-valued $r$-form on $\f E$ with coordinate
expression $\phi = \phi_{\rho\lam_1\dots\lam_r}\,d^\rho\ten(
d^{\lam_1}\wed
\dots \wed d^{\lam_r})$. Then
\beq
    d_K \phi =(\der_\mu \phi_{\rho\lam_1\dots\lam_r}
  + \phi_{\sig\lam_1\dots\lam_r}\, K_\mu{}^\sig{}_\rho)\, d^\rho\ten
  (d^\mu\wed d^{\lam_1}\wed \dots \wed d^{\lam_r})
        \,.
\eeq
\vglue-1.3\baselineskip
\ePf

\bRm\label{Rm1.2}
Now we shall apply $L[K]$ and $d_K$ on specific situation
of the scaled metric $g$.
The metric $g$ can be considered to
be a $\B L^2\ten T^*\f E$-valued
1-form on $\f E$.
Then the covariant exterior
differential $d_K g$ is a $\B L^2\ten T^*\f E$-valued
2-form
defined for any vector fields $X,Y,Z$ by
\beq
        (d_K g)(X,Y)(Z)
=
        \big(\nabla_{X}(Y\Fla)
  -\nabla_{Y}(X\Fla)
  -([X,Y]\Fla)\big)(Z)\,,
\eeq
where ${ }\Fla$ denotes the musical mapping
$g\Fla:T\f E\to\B L^2\ten T^*\f E$.
We have the coordinate expression
\bEq\label{Eq1.9}
        d_K g
= (\der_\lam g_{\rho\mu}
  + g_{\sig\mu}\, K_\lam{}^\sig{}_\rho)
   \, d^\rho \ten( d^\lam\wed d^\mu)\,.
\eEq

On the other hand the musical mapping
$g\Fla$
can be considered as a linear horisontal 1-form
on $T\f E$ with the coordinate expression
$
        g\Fla
=
        g_{\lam\mu} \, \dt x^\lam\, d^\mu\,.
$
Then we have the coordinate equation
\bEq\label{Eq1.10}
    L[K]\, g\Fla = (\der_\lam g_{\rho\mu} \, \dt x^\rho
  + g_{\rho\mu}\, K_\lam{}^\rho)\, d^\lam\wed d^\mu
\eEq
and, if $K$ is linear,
\bEq\label{Eq1.11}
    L[K]\, g\Fla = (\der_\lam g_{\rho\mu}
  + g_{\sig\mu}\, K_\lam{}^\sig{}_\rho)
   \, \dt x^\rho \, d^\lam\wed d^\mu\,,
\eEq
i.e., in  the linear case, $L[K]\, g\Fla$ is a linear horisontal
2-form on $T\f E$
which can be considered to be a $\B L^2\ten T^*\f E$ valued 2-form
on $\f E$ which coincides with $d_K\, g$.
\hfill\ENDE
\eRm

\myssec{Spacetime 2--forms and 2--vectors}
\label{Spacetime 2--forms}

The map
$
T\tau[\f E] : TT\f E \to T\f E \,,
$
can be regarded as a vector valued 1-form
$
\ups : T\f E \to
T^*T\f E \uten{T\f E} T\f E \,,
$
with coordinate expression
$
\ups = d^\lam  \ten \der_\lam \,.
$

We define the {\em spacetime 2--form\/} of
$T\f E$
associated with
$g$
and a spacetime connection
$K $
to be the scaled 2--form
\beq
\Ups[g, K] \byd
g \con \big(\nu[K] \wed \ups\big) :
T\f E \to \B L^2 \ten \Lam^2 T^*T\f E \,.
\eeq
We have the coordinate expression
\bEq\label{Eq1.12}
\Ups[g, K] =
g_{\lam\mu} \, (\dt d^\lam -
K\co\nu\lam \, d^\nu) \wed d^\mu \,
\eEq
and, if $K$ is  linear,
\bEq\label{Eq1.13}
\Ups[g, K] =
g_{\lam\mu} \, (\dt d^\lam -
K\col\nu\lam\rho \, \dt x^\rho \, d^\nu) \wed d^\mu \,
\eEq

We define the {\em spacetime 2--vector\/} of
$T\f E$
associated with
$g$
and a spacetime connection
$K $
to be the scaled 2--vector
\beq
\Lam[g, K] \byd
\ba g \con \big(K \wed \vartheta\big) :
T\f E \to \B L^{*2} \ten \Lam^2 TT\f E \,.
\eeq
We have the coordinate expression
\bEq\label{Eq1.14}
\Lam[g, K] =
g^{\lam\mu} \, (\der_\lam +
K_\lam{}^\nu \, \dt\der_\nu) \wed \dt\der_\mu
\eEq
and, if $K$ is linear,
\bEq\label{Eq1.15}
\Lam[g, K] =
g^{\lam\mu} \, (\der_\lam +
K\col\lam\nu\rho \, \dt x^\rho \, \dt\der_\nu)
\wed \dt\der_\mu \,.
\eEq

\bLm\label{Lm1.3}
We have
\beq
i(\Lam[g,K]) \Ups[g,K] = -4\,.
\eeq
\eLm

\bPf
We have
\beq
i(\Lam[g,K]) \Ups[g,K] = -g^{\lam\mu} \, g_{\lam\mu} = -4\,.
\eeq
\vglue-1.3\baselineskip
\ePf

\mysec{Induced structures on the tangent bundle of the spacetime}
\setcounter{equation}{0}

We study symplectic and Poisson structures
induced on the tangent bundle of the spacetime by the metric $g$
and a spacetime connection $K$.

\myssec{General spacetime connection case}

Let us assume a spacetime connection $K$ given by
(\ref{Eq1.1})
, the spacetime 2--form $\Ups[g,K]$ given by
(\ref{Eq1.12}) and the spacetime 2--vector
$\Lam[g,K]$ given by (\ref{Eq1.14}).

\bLm\label{Lm2.1}
$\Ups[g,K]$ is closed if
and only if the following two conditions are satisfied
\bAl\label{Eq2.1}
        \der_\nu g_{\lam\mu} \,
  + g_{\rho\mu} \, \dt\der_\lam K_\nu{}^\rho
  - \der_\mu g_{\lam\nu} \,
  - g_{\rho\nu} \, \dt\der_\lam K_\mu{}^\rho
& = 0
\\[2mm]
\label{Eq2.2}
        R[K]_{\lam\mu\nu}
  +
        R[K]_{\mu\nu\lam}
  +
        R[K]_{\nu\lam\mu}
& =
        0\,,
\end{align}
where we have set
$R[K]_{\lam\mu\nu} = g_{\rho\nu}\,R[K]_{\lam\mu}{}^\rho$.
\eLm

\bPf
It
follows immediately from the coordinate expression
\bal
        d\Ups[g,K]
& =
  - (\der_\lam g_{\rho\nu}\,K_\mu{}^\rho
  + g_{\rho\nu}\, \der_\lam K_\mu{}^\rho) \,
        d^\lam \wed d^\mu \wed d^\nu
\\
&\quad
  - (\der_\mu g_{\lam\nu} \,
  + g_{\rho\nu} \, \dt\der_\lam K_\mu{}^\rho) \,
        \dt d^\lam \wed d^\mu \wed d^\nu
\\[2mm]
& =
        \fr12\, R[K]_{\lam\mu\nu}\, d^\lam\wed d^\mu\wed d^\nu
\\
&\quad
  - (\der_\mu g_{\rho\nu}
  + g_{\sig\nu}\, \dt\der_\rho K_\mu{}^\sig) \,
        (\dt d^\rho
  - K_\lam{}^\rho \, d^\lam) \wed d^\mu\wed d^\nu\,.
\qquad
\end{align*}
\vglue-1.3\baselineskip
\ePf

Now we shall describe the geometrical interpretation of
the equations (\ref{Eq2.1}) and (\ref{Eq2.2}).
Let us consider the Liouville vector field
$
        I
=
        \dt x^\lam \, \dt\der_\lam\,.
$

\bLm\label{Lm2.2}
    The conditions (\ref{Eq2.1})  and (\ref{Eq2.2})
are equivalent with
\bEq\label{Eq2.3}
    L[I]\, L[K] \, g\Fla = 0\,
\eEq
and
\bEq\label{Eq2.4}
    L[K]\, L[K] \, g\Fla = 0\,,
\eEq
respectively.
\eLm

\bPf
We have
\bal
        L[I]\, L[K] \, g\Fla
& =
  \big(i(I)\,d + d\,i(I)\big)\, \big( (\der_\lam g_{\rho\mu}\, \dt x^\rho
        + g_{\rho\mu}\, K_\lam{}^\rho)\, d^\lam \wed d^\mu
  \big)
\\
& =
        \dt x^\rho \, (\der_{\lam} g_{\rho\mu}
  + g_{\sig\mu}\,\dt\der_\rho K\co\lam\sig)\, d^\lam \wed d^\mu
\,.
\end{align*}
It is easy to see that $L[I]\,L[K]\,g\Fla=0$ if and only if
the condition (\ref{Eq2.1}) is satisfied.

Further from (\ref{Eq1.5}) we have
\beq
    L\big[R[K]\big] \, g\Fla =
        g_{\rho\nu}\, R[K]_{\lam\mu}{}^\rho\,
        d^\lam\wed d^\mu\wed d^\nu
\eeq
i.e., the condition (\ref{Eq2.2}) is equivalent with
$L\big[ R[K]\big] \, g\Fla =0$. But, from (\ref{Eq1.4}),
\bal
L\big[ R[K]\big] \, g\Fla & = - L\big[ [K,K]\big] \, g\Fla
        = - 2\, L[K]\, L[K]\, g\Fla.
\end{align*}
Hence
 (\ref{Eq2.2}) is equivalent with
$L[ K]\,L[K] \, g\Fla =0$.
\ePf


\bLm\label{Lm2.3}
The Schouten bracket
\beq
\big[\Lam[g,K],\Lam[g,K]\big]:
T\f E \to \B L^{*4} \ten \Lam^3 TT\f E
\eeq
has the coordinate expression
\bal
        \big[\Lam[g,K],\Lam[g,K]\big]
& =
        2\, g^{\rho\nu} \, (\der_\rho g^{\lam\mu}
  - g^{\sig\lam} \, \dt\der_\sig K_\rho{}^\mu) \,
        (\der_\lam + K_\lam{}^\kap \, \dt\der_\kap)
        \wed \dt\der_\mu \wed \dt\der_\nu
\\
&\quad
  + R[K]^{\kap\mu\nu}\,
        \dt\der_\kap \wed \dt\der_\mu \wed \dt\der_\nu\,,
\end{align*}
where we have set
$R[K]^{\lam\mu\nu} = g^{\lam\rho}\,g^{\mu\sig}
        \,R[K]_{\rho\sig}{}^\nu$.
\eLm

\bPf
We have
\beq
        i(\big[\Lam[g,K],\Lam[g,K]\big])\,\bet
=
        2\, i(\Lam[g,K]) \, d i(\Lam[g,K])\, \bet
\eeq
for any closed 3-form $\bet$.

Then
\bal
        \big[\Lam[g,K],\Lam[g,K]\big]
& =
        2\, g^{\rho\nu} \, (\der_\rho g^{\lam\mu}
  - g^{\sig\lam} \, \dt\der_\sig K_\rho{}^\mu)
          \, \der_\lam\wed \dt \der_\mu\wed \dt \der_\nu
\\
& \quad +
         \bigg(
        g^{\ome\nu} \, g^{\rho\mu} \, R_{\rho\ome}{}^\kap
  + 2\,g^{\sig\nu}\, K_\rho{}^\kap
  \, (\der_\sig g^{\rho\mu}
  - g^{\rho\ome}\, \dt\der_\ome K_\sig{}^\mu)
        \bigg)\,
        \dt\der_\kap\wed \dt\der_\mu\wed\dt\der_\nu
\\
& =
        2\, g^{\rho\nu} \, (\der_\rho g^{\lam\mu}
  - g^{\sig\lam} \, \dt\der_\sig K_\rho{}^\mu) \,
        (\der_\lam + K_\lam{}^\kap \, \dt\der_\kap)
        \wed \dt\der_\mu \wed \dt\der_\nu
\\
&\quad
  + R[K]^{\kap\mu\nu}\,
        \dt\der_\kap \wed \dt\der_\mu \wed \dt\der_\nu\,.
\end{align*}
\vglue-1.3\baselineskip
\ePf

\bLm\label{Lm2.4}
$\big[\Lam[g,K],\Lam[g,K]\big] = 0$ if and only if the conditions
(\ref{Eq2.1}) and (\ref{Eq2.2}) are satisfied, i.e.,
if and only if the conditions
(\ref{Eq2.3}) and (\ref{Eq2.4}) are satisfied.
\eLm

\bPf
From Proposition \ref{Lm2.3} it follows that
$\big[\Lam[g,K],\Lam[g,K]\big] = 0$ if and only if
\bAl\label{Eq2.6}
        g^{\rho\nu} \, (\der_\rho g^{\lam\mu}
  - g^{\sig\lam} \, \dt\der_\sig K_\rho{}^\mu)
  - g^{\rho\mu} \, (\der_\rho g^{\lam\nu}
  + g^{\sig\lam} \, \dt\der_\sig K_\rho{}^\nu)
& =
        0
\\[1mm]
\label{Eq2.7}
        R[K]^{\kap\mu\nu}
  + R[K]^{\mu\nu\kap}
  + R[K]^{\nu\kap\mu}
& =
        0\,.
\end{align}
But by decreasing the indices in (\ref{Eq2.6})
we get from
$\der_\rho g^{\lam\mu}
=
- g^{\lam\tau}\, g^{\mu\ome} \, \der_\rho g_{\tau\ome}$
just (\ref{Eq2.1}) and by decreasing indices in
(\ref{Eq2.7}) we get just (\ref{Eq2.2})\,.
\ePf


\bTh\label{Th2.5}
The metric $g$ and a general spacetime connection $K$
induce on $T\f E$ natural symplectic and natural Poisson
structures if and only if the conditions (\ref{Eq2.3})
and (\ref{Eq2.4}) are satisfied.
\eTh

\bPf
The regularity of $g$ implies that $\Ups[g,K]$ and $\Lam[g,K]$
are non degenerate. Lemmas
\ref{Lm2.1}, \ref{Lm2.2}
and
\ref{Lm2.4}
then imply that $\Ups[g,K]$ and $\Lam[g,K]$ define
symplectic and Poisson structures, respectively, if and only if
(\ref{Eq2.3})
and (\ref{Eq2.4}) are satisfied.
\ePf

\myssec{Linear spacetime connection case}

We assume a linear spacetime connection $K$.



\bLm\label{Lm2.6}
Let $K$ be a linear spacetime connection.
$\Ups[g,K]$ is closed if
and only if $L[K]\, g\Fla=0$.
\eLm

\bPf
By Lemmas  \ref{Lm2.1} and \ref{Lm2.2}
$\Ups[g,K]$ is closed if and only if (\ref{Eq2.3}) and
(\ref{Eq2.4}) are satisfied. But for a linear spacetime connection
$K$ the horisontal 2-form $L[K]\, g\Fla$ is linear.
Moreover, for any linear horisontal $r$-form $\phi$, we have
$L[I]\phi=\phi$, i.e.,
\beq
   L[I]\, L[K]\, g\Fla = L[K]\, g\Fla\,.
\eeq
Then the condition (\ref{Eq2.3}) is equivalent with
$L[K]\, g\Fla = 0$ which implies
(\ref{Eq2.4}).
\ePf

\bRm   \label{Rm2.7}
In \cite{Jan95} we have proved that the spacetime 2-form
$\Ups[g,K]$ is closed if and only if $d_K g=0$. By Remark \ref{Rm1.2}
it coincides with Lemma  \ref{Lm2.6}.
\hfill\ENDE
\eRm

\bLm\label{Lm2.8}
Let $K$ be a linear spacetime connection then
$\big[\Lam[g,K],\Lam[g,K]\big] = 0$ if and only if
$L[K]\, g\Fla =0.$
\eLm

\bPf
It follows immediately.
\ePf


\bTh\label{Th2.11}
Let $K$ be a linear spacetime connection.
Then the following identities are equivalent:

\medskip

(1) $L[K]\, g\Fla =0\,$.

\medskip
(2) $d_K g = 0\,.$

\medskip
(3) $d\Ups[g,K] = 0\,.$

\medskip
(4) $\big[\Lam[g,K],\Lam[g,K]\big] = 0\,.$
\eTh

\bPf
    It follows from Remark \ref{Rm1.2} and Lemmas \ref{Lm2.6}
and  \ref{Lm2.8}.
\ePf

\bCr \label{Cr2.12}
A linear spacetime connection $K$ and the metric $g$ induce
on $T\f E$ natural symplectic and Poisson structures
if and only if $d_K g = 0 = L[K]\, g\Fla\,.$
\hfill\ENDE
\eCr

\bLm\label{Lm 2.9}
Let $K$ be a linear spacetime connection then
the following three identities are equivalent:

(1) $d_K g = 0$.

(2) $L[K]\ g\Fla =0$.

(3) $(\nab_X g)(Y,Z) - (\nab_Y g)(X,Z)
=
  2\, g(\tau[K](X,Y),Z)\,.$
\eLm

\bPf
(1) $\Leftrightarrow$ (2).
It follows immediately from Remark \ref{Rm1.2}.

(1) $\Leftrightarrow$ (3). Let us recall that
for a linear connection  $K$ we have
\bEq
        2\, \tau[K](X,Y)
=
        \nab_Y X
  - \nab_X Y
  + [X,Y]\,.
\eEq
Then, by Remark \ref{Rm1.2},
\bal
        (d_K g)(X,Y)(Z)
& =
 (\nab_X g)(Y,Z)
  + g(\nab_X Y, Z)
  - (\nab_Y g)(X,Z)
  - g(\nab_Y X, Z)
  - g([X, Y], Z)
\\
& = (\nab_X g)(Y,Z)
  - (\nab_Y g)(X,Z)
  + g(\nab_X Y - \nab_Y X - [X, Y], Z)
\\
& = (\nab_X g)(Y,Z)
  - (\nab_Y g)(X,Z)
  - 2\, g(\tau[K](X, Y), Z)
\, .
\end{align*}
\vglue-1.3\baselineskip
\ePf

\bCr\label{Cr 2.10}
If $K$ is a torsion free connection then $d_K g =0 = L[K] g\Fla$
is equivalent to
$(\nab_X g)(Y,Z)
  = (\nab_Y g)(X,Z)$, i.e., for a torsion free linear
connection the (0,3)-tensor
field  $\nab g$ is symmetric.
\hfill\ENDE
\eCr

\bTh \label{Th2.13}
Let $K$ be a linear torsion free spacetime connection.
Then the following identities are equivalent:

\medskip
(1)
$
        \nab g
$             is a symmetric (0,3)-tensor field.

\medskip
(2) $d\Ups[g,K] = 0\,.$

\medskip
(3) $\big[\Lam[g,K],\Lam[g,K]\big] = 0\,.$
\eTh

\bPf
By Corollary \ref{Cr 2.10} for a linear torsion free
connection the identity
$d_K g = 0 = L[K]\, g\Fla$ is equivalent  with
$\nab g$ to be fully symmetric.
\ePf

\bCr\label{Cr2.14}
A linear torsion free spacetime connection $K$
and the metric $g$ induce on $T\f E$
natural symplectic and Poisson structures
if and only if the covariant differential
$\nab g$ is a symmetric (0,3)-tensor field.
\hfill\ENDE
\eCr

\bRm \label{Rm2.15}
Let us assume the torsion free space time metric connection $K[g]$
given by the Cristtoffel symbols (\ref{Eq1.3}).
Then $\nabla g = 0$ (i.e. symmetric)
and we have the canonical natural symplectic
and Poisson structures on $T\f E$ given by
$\Ups[g]=\Ups\big[g,K[g]\big]$
and
$\Lam[g] = \Lam\big[g,K[g]\big]$.
Moreover, in the metric case,
$\Ups[g]= dg\Fla.$
\hfill\ENDE
\eRm


\bibliographystyle{amsplain}

\begin{thebibliography}{10}
\bibitem{Jan95}
        J. Jany\v ska,
         {\em Remarks on symplectic and contact 2--forms in
         relativistic theories,}
         Bollettino U.M.I. (7) 9--B (1995), 587--616.
\bibitem{Jan96}
        J. Jany\v ska,
        {\em Natural symplectic structures on the tangent
                bundle of a space--time,}
        The Proceedings of the Winter School
                Geometry and Topology (Srn\'\i , 1995),
                Supplemento ai Rendiconti del Circolo Matematico
                di Palermo, Serie II -- No. 43 (1996),
         153--162.
\bibitem{JanMod96a}
        J. Jany\v ska, M. Modugno,
        {\em Classical particle phase space in general
                relativity,}
        Differential Geometry and Applications,
                Proc. Conf., Aug. 28 -- Sept. 1, 1995, Brno,
                Czech republic, Masaryk University, Brno 1996,
                573--602.
\bibitem{JanMod97}
        J. Jany\v{s}ka, M. Modugno,
        {\em On quantum vector fields in general relativistic quantum
        mechanics,}
        in: Proc. 3rd Internat. Workshop Diff.
        Geom. and its Appl., Sibiu (Romania) 1997,
        General Mathematics {\bf 5} (1997), 199--217.
\bibitem{Jan00}
        J. Jany\v ska,
         {\em Natural Poisson and Jacobi structures on the tangent
        bundle of a pseudo-Riemannian manifold,}
        Contemporary Mathematics {\bf 288}, Global
        Diff. Geom.: The Math. Legacy of Alfred Gray,
        eds. M. Fern\'andes
        and J. A. Wolf, 343--347.
\bibitem{Jan01}
        J. Jany\v ska,
        {\em Natural vector fields and 2-vector fields
           on the tangent bundle of a pseudo-Riemannian manifold,}
        Archivum Mathematicum (Brno) {\bf 37} (2001),
        143--160.
\bibitem{KruJan90}
        D.~Krupka, J.~Jany\v{s}ka,
        { Lectures on Differential
        Invariants}, Folia Fac. Sci. Nat. Univ. Purkynianae Brunensis,
        Brno, 1990.
\bibitem{KolMicSlo93}
        I. Kol\' a\v r, P. W. Michor, J. Slov\' ak,
         {Natural Operations in  Differential Geometry,}
         Springer--Verlag 1993.
\bibitem{LibMar87} P. Libermann, Ch. M. Marle,
        {Symplectic Geometry and Analytical Mechanics,}
        Reidel Publ., Dordrecht 1987.
\bibitem{Nij72}
        A.~Nijenhuis,
        {\em Natural bundles and their general properties,}
        Diff. Geom., in honour of K. Yano, Kinokuniya, Tokyo 1972,
        317--334.
\bibitem{Vai94} I.~Vaisman,
        {Lectures on the Geometry of Poisson Manifolds},
        Birkh\"auser, Verlag 1994.

\end{thebibliography}
\makeatletter \renewcommand{\@biblabel}[1]{\hfill#1.}\makeatother

\end{document}